\documentstyle[amsmath,amsbsy,amsfonts]{article}
\numberwithin{equation}{section}

\newtheorem{theorem}{Theorem}[section]  
\newtheorem{definition}[theorem]{Definition}

\newtheorem{lemma}[theorem]{Lemma}  
\newtheorem{proposition}[theorem]{Proposition}  
\newtheorem{corollary}[theorem]{Corollary}  
\newtheorem{remark}[theorem]{Remark}

\newcommand{\A}{{\cal A}} 
\newcommand{\B}{{\cal B}}

\newcommand{\V}{{\cal V}}
\newcommand{\Ss}{{\cal S}}

\newcommand{\depth}{\mbox{\rm depth}}
\newcommand{\rk}{\mbox{\rm rk}}

\newcommand{\proof}{{\bf Proof.~}}
\newcommand{\qed}{~~\mbox{$\Box$}}

\begin{document}

\title{Taylor and minimal resolutions of homogeneous polynomial ideals}
\author {
{\sc Sergey Yuzvinsky }\\
{\small\it University of Oregon,
Eugene, OR 97403 USA}\\ 
{\small\it yuz@math.uoregon.edu}}

\date{May 14, 1999}
\maketitle

\section{Introduction}
\bigskip
In the theory of monomial ideals of a polynomial ring $S$ over a field $k$,
 it is convenient that
for each such ideal $I$ there is a standard free resolution, so called Taylor
resolution, that can be canonically constructed from the minimal system of
monomial generators of $I$ (see \cite{Es}, p.439 and section 2). 
On the other hand no construction of a minimal resolution for an arbitrary
monomial ideal has been known. Recently a minimal resolution was constructed in
\cite{BPS} for a class of so called generic monomial ideals. Also in
\cite{BC,GPW1,GPW2} various invariants of monomial ideals were related to
combinatorics of the lattice $D$ of the least common multiples (lcm) of generating
monomials. In particular in \cite{GPW2}
the Betti numbers of the $S$-module $S/I$ were expressed through homology of $D$
and it was proved that even the algebra structure of 
${\rm Tor}_*^S(S/I,k)$ 
was defined by that lattice although explicit formula was
not given in that paper.

Given a system of generators of an arbitrary ideal $I$ of $S$, one can 
factor the generators in irreducibles and construct
the Taylor complex similarly to the Taylor resolution of a monomial ideals. In
general this complex is not acyclic. One non-monomial case where it is acyclic
was used in \cite{Yu1}.

In the present paper, a necessary and sufficient condition is given
for the Taylor complex
of a system $\A$ of homogeneous polynomials to be acyclic (Theorem \ref{taylor}). 
This condition involves
the local homology of the lattice $D$ of lcm of elements form $\A$
and the depth of ideals generated by their irreducible factors. 
If this condition holds then the Betti numbers of $S/I$ are
defined by the local homology of $D$
 similarly to the case of monomial ideals (Theorem
\ref{betti}).
Moreover in section 3
we exhibit a DGA defined by combinatorics whose cohomology algebra is
isomorphic to the algebra 
${\rm Tor}_*^S(S/I,k)$. This construction makes sense for arbitrary graded
lattice (see definition in section 3) and generalizes the DGA constructed in
\cite{Yu2, Yu3} for cohomology algebra of complex subspace complement.
Section 4 contains that can be considered the main result of the
paper (Theorem \ref{minimal}). There for any ideal $I$ having the Taylor
resolution we give a combinatorial construction of
a subcomplex of it that is a minimal resolution of $S/I$.
This construction is not canonical and involves computations of homology of 
posets which hardly can be avoided in general. We also describe completely
the class of ideals for which 
our minimal resolution reduces to the minimal resolution from \cite{BPS}.
Finally in section 5 we give examples of classes of $\A$
satisfying the condition of Theorem \ref{taylor}.
For instance we consider the ideals whose generators are
products of linear polynomials
which is important for theory of hyperplane and subspace arrangements.

The author is grateful to D.Eisenbud for discussions of the results of the paper.
In particular the proof of Theorem \ref{minimal} would have been longer without
his advice.
\bigskip
\section{Taylor complex}
\bigskip
 Let $S=k[x_1,\ldots,x_{n}]$ be the polynomial ring
over a field $k$,
$\A=\{Q_1,\ldots,Q_m\}$ a set of 
homogeneous polynomials from $S$, and $I=I(\A)$ 
the ideal of $S$ generated by $\A$. 
In this section we give a necessary and sufficient condition for the Taylor complex
of $\A$ to be a resolution of the $S$-module
$S/I$.
We can assume that $\A$ generates $I$ minimally, in particular none of
polynomials from $\A$
divides another one. Also we will always regard all the subsets of $\A$
provided with the ordering induces by a fixed linear ordering on $\A$. 

Let $K: 0\to K_{m}\to \cdots \to K_1\to K_{0}=k\to 0$ be the
augmented shifted by -1 chain complex over $k$ of the simplex whose set of vertices
is $\A$.
Explicitly the linear space $K_p$ has a basis consisting of all the 
subsets of $\A$ with $p$
elements. The differential $d_p:K_p\to K_{p-1}$ is given with respect to
these bases by the
matrix with entries $d_{\sigma,\tau}$ ($\sigma,\tau\in \A, |\sigma|=p-1, |\tau|=p$)
equal to $(-1)^{\epsilon(\tau,i)}$ if $\sigma=\tau\setminus\{Q_i\}$ where
$\epsilon(\tau,i)=|\{Q_j\in\tau|j<i\}|$ and 0 if $\sigma\not\subset\tau$.
For each $\sigma\subset\A$ denote by $Q_{\sigma}$ the lcm
of all $Q_i\in\sigma$.

\begin{definition}
The Taylor complex of $\A$ is the complex $\tilde K$ of the free $S$-modules $\tilde
K_p=K_p\otimes S$ with the differentials $\tilde d_p:\tilde K_p\to\tilde K_{p-1}$
given by the matrix $\tilde d_{\sigma,\tau}=(Q_{\tau}/Q_{\sigma})d_{\sigma,\tau}$.
\end{definition}

Notice that $H_0(\tilde K)=S/I$.
Our goal is to give equivalent conditions for the acyclisity of $\tilde K$.

We will be dealing with finite
 posets and need to introduce some notation. To every poset
one can assign its complex of flags (i.e., increasing sequences of elements) 
and attribute
the topological invariants of this complex to the poset itself. In this sense we
will use homology of a poset and homotopy equivalence of two of them. If $L$ is
a lattice with the minimal element $\hat 0$ and the maximal one $\hat 1$ then the
flag complex of $L\setminus\{\hat 0,\hat 1\}$ is homotopy equivalent to its
atomic and coatomic complexes (e.g., see \cite{Bj}). 
For instance, recall that the atomic complex is
the abstract complex whose vertices are all the atoms of $L$ and a set of atoms
is a simplex if it is bounded above in $L\setminus\{\hat 1\}$.
If $\sigma\subset L$ then $\bigvee(\sigma)$ denotes the least upper
bound (join) of $\sigma$. 
For every $X\in L$ we put $L_{\leq X}=(\hat 0,X]=\{Y\in L|0<Y\leq X\}$ and
$L_{<X}=L_{\leq X}\setminus\{X\}$. 

Now denote by $D=D(\A)$ the lattice of the lcm of all 
subsets of $\A$ ordered by divisibility.
This lattice is provided with the monotone map $\phi: \B\to D$ from the Boolean
lattice $\B$ of all the subsets of $\A$ given by $\phi(\sigma)=
Q_{\sigma}$.
Denote by $D_0$ the poset $D$ with the smallest element (constant polynomial 1)
deleted and call a subset of a poset {\it decreasing} if with each element it
contains all smaller ones.

\begin{lemma}
\label{homotopic}
For every decreasing subset $F$ of $D_0$ 
the restriction of
 $\phi$ to $\phi^{-1}(F)$ is a homotopy equivalence. 
\end{lemma}
\proof
For $Q\in F$ put $\B(\leq Q)=\phi^{-1}(\{P\in F|P\leq Q\})$ and notice that
the poset $\B(\leq Q)$ has the unique maximal element equals $\{Q_i\in\A|Q_i {\rm \
divides}\ Q\}$. Thus $\B(Q)$ is contractable and the restriction of
$\phi$ is a homotopy
equivalence (cf. \cite{Qu}).
\qed

\medskip
The lattice $D$ is naturally embedded as a join-sublattice in a larger polynomial
 lattice $W$. Let $E$ be the 
set of all irreducible factors of elements form $\A$ (chosen one from each
equivalence class) and $W$ the lattice of all the products of
elements from $E$ ordered by divisibility (repetitions are allowed).
 Clearly $W$ is isomorphic as a partially
ordered semigroup to ${\mathbb{N}}^k$ ($k=|E|$).

To study the complex $\tilde K$ we use the evaluation at vectors $v\in\bar k^n$.
This means that for every $v\in\bar k^n$ we consider the complex $K(v)$ with
$K(v)_p=K_p\otimes \bar k$ and $d(v)_p:K(v)_{p+1}\to K(v)_p$ given by the matrix
$d(v)_{\sigma,\tau}=\tilde d_{\sigma,\tau}(v)$.
On the other hand, every $v$ defines a subset $E_v$ of $E$ of all the elements of
$E$ that vanish at $v$. Denote by $W(v)$ the sublattice of $W$ of all the elements
of $W$ whose all irreducible factors are from $E_v$ and by $\psi_v$ the monotone
map (``projection'') $D\to W(v)$ defined by
$$\psi_v(Q)=\max\{P\in W(v)|P\leq Q\}.$$
For each $P\in W(v)$ we put $D(v,<P)=\{Q\in D|\phi_v(Q)<P\}$.
Besides every $P\in W(v)$ defines the
subcomplex $K_{v,P}$ of $K\otimes \bar k$ spanned by $\sigma\subset\A$ such that
$\psi_v(Q_{\sigma})=P$.
Here and in the rest of the paper we consider only $P\in W(v)$ such that $P\in
\psi_v(D)$. Substituting in the last definition ``='' by ``$\leq$'' or
``$<$'' we obtain complexes $K_{v,\leq P}$ and $K_{v,<P}$ respectively. 
Clearly $K_{v,P}=K_{v,\leq P}/K_{v,<P}$.

Notice that if $E(v)=\emptyset$ then $W(v)=\{1\}$ and $K(v)\simeq K$.
We need the following lemma.

\begin{lemma}
\label{eval}
(i) $K(v)\simeq\oplus_{P\in \psi_v(D)}K_{v,P}$ for every $v\in\bar k^n$.

(ii) The complex $K_{v,\leq P}$ is acyclic for every $v$ and $P\in \psi_v(D)$.

(iii) The complex $K_{v,<P}$ is homotopy equivalent to the chain complex of 
$D(v,<P)$
shifted by -1.
\end{lemma}

\proof (i) The decomposition of $K(v)$ in the direct sum according to 
$\psi_v(Q_{\sigma})$ follows from definitions.
The required isomorphism can
 be given by $\sigma\mapsto (Q_{\sigma}/P)(v)\sigma$ where
$P=\psi_v(Q_{\sigma})$.

(ii) The complex $K_{v,\leq P}$ is the chain complex of the simplex on the vertices
$Q_i\in\A$ with $\psi_v(Q_i)\leq P$. The result follows.

(iii) The complex $K_{v,<P}$ is the shifted by -1 chain complex of an abstract complex
whence it is homotopy equivalent to the poset of its non-empty simplexes ordered by
inclusion. More explicitly $K_{v,<P}$ is homotopy equivalent (after the shift)
to the subposet of $\B$
equal to $\phi^{-1}(D(v,<P))$. Since $D(v,<P)$ is decreasing we can apply 
Lemma \ref{homotopic} which completes the
proof.                 \qed

\medskip
It will be convenient for us to express properties of $v$ in terms of the set $E_v$
of irreducible polynomials. Notice that these sets can be characterized
intrinsically as subsets $G\subset E$ such that the ideal $J(G)$ generated by $G$
does not contain a power of an element from $E\setminus G$. We will call
those sets {\it
saturated}. For every two $v,v'$
such that $E_v=E_{v'}$ we have $\psi_v=\psi_{v'}$. Thus we will often use $\psi_G$,
$W(G)$,
and $D(G,<P)$ for $P\in W(G)$ instead of $\psi_v$,
$W(v)$, and $D(v,<P)$ respectively for $v$ such
that $E_v=G$.
\begin{theorem}
\label{taylor}
The Taylor complex of $\A$ is a resolution of $S/I(\A)$ if and only if 
for every saturated set $G\subset E$ and $P\in \psi_G(D)$
we have $\tilde H_p(D(G,<P)))=0$ for every $p$ such that $p\geq\depth(J(G))-1$.
\end{theorem}

\proof
According to the celebrated Buchsbaum - Eisenbud criterion
(see for example \cite{No}, Sect. 6.4, Theorem 15)
 the complex $\tilde K$ is a resolution
of $S/I$ if and only if the following three conditions hold:

(a)$\depth F(\tilde d_p)\geq p$ for every $p=1,2,\ldots$ where $F(\tilde d_p)$ (or
shortly $F_p$) is
the Fitting ideal of $\tilde d_p$;

(b)$\rk\tilde d_m=\rk\tilde K_m$;

(c)$\rk\tilde d_{p+1}+\rk\tilde d_p=\rk \tilde K_p$ for $1\leq p<m$.

Consider $K(v)$ for a $v\in\bar k^n$. Clearly we have $\rk \tilde K_p=\dim K_p$
and $\rk d_p(v)\leq\rk_S\tilde d_p$ for all $p$.
Using these and choosing $v$ so that $E_v=\emptyset$ we 
can obtain by induction on $p$ that 
$\rk_S\tilde d_p=\rk d_p$ and the conditions (b) and (c)
(cf. \cite{Yu1}, Theorem 1.3).
Thus we will 
focus on the condition (a). 

For an arbitrary $v$,
Lemma \ref{eval} (i) and (ii) 
and the exact sequence of the pair $(K_{v,\leq P},K_{v,<P})$ ($P\in \psi_v(D)$)
give for all $p$
\begin{equation}
\label{2.2}
H_p(K(v))=\oplus_{P\in \psi_v(D)}H_{p-1}(K_{v,<P}).
\end{equation}

Then Lemma \ref{eval} (iii) gives the
following equivalent form of \eqref{2.2}
\begin{equation}
\label{2.3}
H_p(K(v))=\oplus_{P\in \psi_v(D)}\tilde H_{p-2}(D(v,<P)).
\end{equation}

Now we include the Fitting ideals into consideration. Clearly the statement
$H_p(K(v))\not=0$ is equivalent to $\rk d_i(v)<\rk d_i$ for $i=p$ or $i=p+1$ which
is equivalent to the inclusion $v\in\V(F_{p+1})\cup\V(F_p)$. Here and in the rest of
the paper for any ideal $J$
of $S$ we denote by $\V(J)$ its variety in $\bar k^n$. In particular the left hand
side of \eqref{2.3} is nonzero for $v$ from a closed algebraic
subset of
$\bar k^n$. The right hand side of \eqref{2.3} is nonzero for $v$ if there exists
$P\in \psi_v(D)$ (whence $v\in \V(J(P))$) such that 
\begin{equation}
\label{2.3'}
\tilde H_{p-2}(D(v,<P))\not=0.
\end{equation}
If the condition \eqref{2.3'} holds for one $v_0\in J(\V(P))$ it may happen that it
holds only for $v$ from a closed algebraic subset of smaller dimension. On the
other hand if we put $G=E_{v_0}$ then \eqref{2.3'} holds for the same $P$ and $v$ from
the open dense subset of $\V(J(G))$ (of all $v$ such that $E_{v}=G$).
Thus \eqref{2.3} implies
\begin{equation}
\label{2.4}
\V(F_{p+1})\cup\V(F_p)=\bigcup_{G\in\Ss_p}\V(J(G))
\end{equation}
where $\Ss_p$ consists of all saturated subsets $G$ of $E$ such that for each of
them there exists $P\in \psi_G(D)$ satisfying \eqref{2.3'} 
for some (whence every) $v$ with $E_v=G$.
 Applying the Nullstellensatz we can
rewrite \eqref{2.4} in the form
\begin{equation}
\label{2.5}
\bar F_{p+1}\cap\bar F_p=\bigcap_{G\in\Ss(p)}\bar J(G)
\end{equation}
where the bar means the radical of an ideal.

Now we are ready to prove 
that the condition (a) above is equivalent to the condition of the theorem.
We are going to use that for any proper ideal there exists containing it prime
ideal of the same depth (see for example \cite{No}, 5.5, Theorem 16).

Suppose that $\depth F_p\geq p$ for every $p$.  Then due to \eqref{2.5}
the same is true for any
prime ideal containing $\bar J(G)$ for $G\in\Ss(p)$ whence for any ideal $J(G)$.
In other words if for some $G\subset E$ as in the condition of the theorem
$r\geq\depth(J(G))-1$ then $G\not\in \Ss(r+2)$, i.e., for every $v$ with
$E_v=G$ and $P\in\psi_G(D)$ we have $\tilde H_r(D(v,<P))=0$. This is precisely the
condition of the theorem.

Conversely suppose that the condition of the theorem holds and consider $G\in\Ss(p)$
for some $p$. Then the assumption implies that $\depth J(G)\geq p$ whence the
same is true for any prime ideal containing either side of \eqref{2.5}. 
This imeddiately implies condition (a) for $p$.                  \qed

\medskip
If a set $\A$ of homogeneous polynomials satisfies the condition of Theorem
\ref{taylor} then it is easy to give a combinatorial interpretation of the Betti
numbers of $S/I(\A)$, i.e., $b_p=\dim {\rm Tor}_p^S(S/I(\A),k)$.

\begin{theorem}
\label{betti}
 Suppose a set $\A$ of homogeneous polynomials from $S$ satisfies the condition of
Theorem \ref{taylor}. Then
$$b_p=\sum_{Q\in D}\dim \tilde H_{p-2}(D_{<Q}).$$
\end{theorem}
\proof
The complex whose homology is ${\rm Tor}_*^S(S/I,k)$, is $\tilde K\otimes_S k$ that
coincides with $K(0)$. Clearly $W(0)=W$ whence $\psi_0$ is the embedding $D\subset
W$. Thus $D(0,<P)=D_{<P}$ for any $P\in D=\psi_0(D)$.
Applying Lemma \ref{eval} for $v=0$ we obtain the result.
                     \qed

\medskip
\begin{remark}
Using Lemma \ref{homotopic} and the exact sequence of a pair we can rewrite the
formula for the Betti numbers as
$$b_p=\sum_{Q\in D}\dim H_{p-1}(K(Q))$$
where the complex $K(Q)$ is generated by
 $\B(Q)=\{\sigma\in\B|Q_{\sigma}=Q\}$ and its differential is the restriction of
$d$. This should be compared with 
Theorem \ref{minimal}.
\end{remark}
\bigskip
\section{Multiplication}
\medskip
First in this section we consider a pure combinatorial set up of a 
grade lattice (see definition below)
and define 
a DGA on the relative atomic complex of this lattice.
A special case of this definition was used in \cite{Yu2}
in order to describe the rational cohomology ring of a subspace complement.
In certain special cases it is known to give 
even the integer cohomology \cite{Fe,ML}.

Let $L$ be a lattice with the minimal element $\hat 0$. The atoms of $L$ are
provided with an arbitrary but fixed linear ordering.
 A {\it grading} of $L$ 
is a strictly monotone map $\rk:L\to {\mathbb{N}}^s$ (for some positive integer $s$)
 with $\rk(0)=0$ and 
$$\rk (X\vee Y)+\rk (X\wedge Y)\leq \rk(X)+\rk(Y)$$
for $X,Y\in L$.
 We call $\rk X$ {\it rank of $X$} for $X\in L$.

Let us recall the relative atomic complex $\Delta=\Delta(L)$ of $L$ (cf. \cite
{Yu3}).
It is a chain complex (over a field or $\mathbb{Z}$) of linear spaces (free modules
resp.) whose $p$-th term $\Delta_p$ has a
basis consisting of subsets $\sigma$ of atoms of $L$ with $|\sigma|=p$ (notice
the unusual grading) and with differential
defined by
$$d\sigma=\sum_{Z_i\in \sigma,\bigvee(\sigma_i)=\bigvee(\sigma)}
(-1)^{\epsilon(\sigma,i)}\sigma_i.$$
(Here as in section 2, $\epsilon(\omega,i)=|\{Z_j\in E|j<i|$ and if
$Z_i\in\omega$ then $\omega_i=\omega\setminus\{Z_i\}$). As usual $\Delta_0$ is
spanned by the empty set of atoms of $L$.
It is easy to see that $\Delta=\oplus_{X\in L}\Delta_X$ where $\Delta_X$ is
spanned by $\sigma$ with $\bigvee(\sigma)=X$ and $\tilde H_p(\Delta_X)=
\tilde H_{p-2}(L_{<X})$ where the latter natural isomorphism is given by the
boundary map in an exact sequence of $(L_{\leq X},L_{<X})$ (with a shift of
dimension).

\begin{proposition}
\label{dga}
The complex $\Delta$ gets the structure of a DGA via the bilinear
multiplication defined by
\begin{equation}
\label{4.1}
\sigma\cdot\tau=\begin{cases}
0& \text{if $\rk(\bigvee(\sigma\cup\tau))\not=\rk(\bigvee(\sigma))
+\rk(\bigvee(\tau))$},\\
(-1)^{\epsilon(\sigma,\tau)}\sigma\cup\tau & \text{otherwise}.
\end{cases}
\end{equation}
where $\epsilon(\sigma,\tau)$ is the
parity of the permutation of $\sigma\cup\tau$
(shuffle) putting all elements of $\tau$ after elements of
$\sigma$ and preserving fixed orders inside these sets.
\end{proposition}
\proof It suffices to check the Leibniz
rule. We consider two cases.

(i) $\rk(\bigvee(\sigma\cup\tau))\not=\rk(\bigvee(\sigma))
+\rk(\bigvee(\tau)).$
Then $\sigma\cdot\tau=0$. On the other hand, for every $\sigma_i$
such that
$\bigvee(\sigma_i)=\bigvee(\sigma)$ 
\begin{eqnarray}
\rk(\bigvee(\sigma_i\cup\tau))\leq\rk\bigvee((\sigma\cup\tau))<\rk(\bigvee
(\sigma))+\rk(\bigvee(\tau))\notag\\
=\rk(\bigvee(\sigma_i))+\rk(\bigvee(\tau))
\end{eqnarray}
whence $d(\sigma)\cdot\tau=0$. Similarly $\sigma\cdot d(\tau)=0$ that gives the
Leibniz rule in this case.

(ii) $\rk(\bigvee(\sigma\cup\tau))=\rk(\bigvee(\sigma))
+\rk(\bigvee(\tau)).$ First notice that in this case $\sigma\cap\tau=\emptyset$.
Indeed $Z\in\sigma\cap\tau$ would imply
$$\rk(\bigvee(\sigma\cup\tau))\leq\rk(\bigvee
(\sigma_i))+\rk(\bigvee(\tau))-\rk Z<\rk(\bigvee(\sigma))+\rk(\bigvee(\tau)).$$

 Further both sides of the Leibniz equality are
combinations of basic elements $\omega_i=(\sigma\cup\tau)_i$. 
Suppose $Z_i\in\sigma$.
Then the coefficient of $\omega_i$ in the left hand side is
\begin{equation}
(-1)^{\epsilon(\omega,i)+\epsilon(\sigma,\tau)}
\end{equation}
and that in the right hand side side is
\begin{equation}
(-1)^{\epsilon(\sigma,i)+\epsilon(\sigma_i,\tau)}.
\end{equation}
We have
\begin{eqnarray}
\epsilon(\sigma,i)+\epsilon(\sigma_i,\tau)\equiv
\epsilon(\sigma,i)+\epsilon(\tau,i)+\epsilon(\sigma,\tau)\notag\\
\equiv \epsilon(\omega,i)
+\epsilon(\sigma,\tau) ({\rm mod}\ 2).
\end{eqnarray}
 which implies the Leibniz equality in this
case. The case where $Z_i\in\tau$ can be completed similarly.
                   \qed

\medskip
Let us consider several important particular cases of this construction.

{\bf Examples.} 1. If $L$ is a geometric lattice graded by its standard rank then
 the complex $\Delta$ over $\mathbb{Z}$
is homotopy equivalent to the
shifted by -1 Whitney complex (e.g., see \cite{OT}, p.142). 
The homology ring of the DGA in this case
is the Orlik-Solomon algebra of $L$ after the shift.

2. If $L$ is the intersection lattice of a complex subspace arrangement graded by
the codimensions of its elements then 
the DGA $\Delta$ is the DGA from \cite{Yu2} whose homology ring is isomorphic with
properly regraded cohomology ring of the subspace complement.

3. Suppose $L=D$ that is the least common multiple lattice from above corresponding
to a set $\A$ of homogeneous polynomials. As in the previous sections denote the
set of atoms of $D$ by $E$ and put $r=|E|$. Fix a linear order on $E$.
 Then one defines the grading $ \rk:D\to {\mathbb{N}}^r$ assigning to
each $Q\in D$ the vector of multiplicities of the irreducible factors of $Q$. 
The complex $\Delta(D)$ is homotopy equivalent to 
the complex $\oplus_{Q\in D}D_{<Q}$ (cf. \cite{Yu3}).

One can generalize this definition using more general ordered semigroups instead
of ${\mathbb{N}}^s$ but we will not use this in the paper.

To show that the example 3 above can be used for computation of the algebra
structure on ${\rm Tor}_*^S(S/I,k)$ one can use the well-known
graded algebra structure on the Taylor complex (see for example \cite{Av}, p. 6).
 To define a (bilinear over $S$) product on $\tilde K$ it is enough to define it on
the standard generators. For $\sigma,\tau\subset\A$ we put $\sigma\cdot\tau=0$ if
$\sigma\cap\tau\not=\emptyset$ and otherwise
\begin{equation}
\label{2.6}
\sigma\cdot\tau=(-1)^{\epsilon(\sigma,\tau)}[Q_{\sigma},
Q_{\tau}]\sigma\cup\tau
\end{equation}
where $[\ ,\ ]$ is the greatest common divisor. 

\begin{theorem}
\label{pairing/complex}
There is a canonical $k$-algebra isomorphism
${\rm Tor}^S_*(S/I,k)\to H_*(\Delta(D))$ where the algebra structure
in the right hand side is induced by the 
DGA structure on $\Delta(D)$ defined in Example 3 above.
\end{theorem} 
\proof 
It is obvious from definitions that $K(0)$ can be
identified with the relative atomic complex $\Delta(D)$.
Since any DGA structure on a free resolution
can be used for defining the multiplicative structure on Tor, we can use the DGA
defined above on $\tilde K$. It is straightforward to check that it induces on 
$\Delta(D)$ the DGA structure defined in Example 3 above.                  \qed

\bigskip
\section{Minimal resolution}
\bigskip
In this section we fix a set $\A$ of homogeneous polynomials generating ideal
$I$ satisfying the condition of Theorem \ref{taylor} and starting from the Taylor
resolution $\tilde K$ of $S/I$ construct its minimal resolution $\tilde M$.
The complex $\tilde M$ will be realized as a subcomplex of $\tilde K$. Moreover
we first find a needed subcomplex $M$ of $K$ and then pass to $\tilde M$ in the
same way as we obtained $\tilde K$ from $K$.

Recall that the acyclic complex $K$ has the standard basis consisting 
of the elements $\sigma$ of the Boolean lattice $\B$ (i.e., subsets of $\A$)
graded by $|\sigma|$. The monotone map $\phi:\B\to D$ defines the partition of
$\B$ into sets $\B(Q)=\phi^{-1}(Q)$ ($Q\in D$).
Denote by $K(Q)$ the graded subspace of $K$ generated by $\B(Q)$ and
notice that $K=\oplus_{Q\in D}K(Q)$ as a graded linear space. If we provide $K(Q)$
with the restriction $d(Q)$ of $d$ it becomes a chain complex. More precisely
$d(Q)(a)=\pi_Qd(a)$ for every $a\in K(Q)$ where $\pi_Q:K\to K(Q)$ is the canonical
projection.

Now we need to make a noncanonical choice. Let $Z(Q)_p$ and $B(Q)_p$ be the
spaces of cycles and boundaries respectively of degree $p$ in $K(Q)$. For each
$p$ there are two exact sequences
$$0\to Z(Q)_p\to K(Q)_p\to B(Q)_{p-1}\to 0$$
and 
$$0\to B(Q)_p\to Z(Q)_p\to H(Q)_p\to 0.$$
We fix a splitting of each of the sequences. In other words we represent
$$K(Q)_p=B(Q)_p\oplus H'(Q)_p\oplus B'(Q)_{p-1}$$
where $B(Q)_p\oplus H'(Q)_p=Z(Q)_p$
with the restriction of $d(Q)$ giving an isomorphism $B'(Q)_{p-1}\to B(Q)_{p-1}$
and the restriction of the projection $Z(Q)_p\to H(Q)_p$ giving an isomorphism 
$H'(Q)_p\to H(Q)_p$.

In the case where the graded
 subspace $H'=\oplus_QH'(Q)$ of $K$ is invariant under $d$
we can take this subspace for $M$. However it is easy to find examples where this
is false (cf. Example 4.2).
 Our goal is to find a subcomplex $M$ of $K$ that is the graph of a degree
-1 linear
map $f:H'\to B'=\oplus_QB'(Q)$.
The map $f$ is defined by the following lemma.
\begin{lemma}
\label{map}
For each 
$a\in H'_p$ there exists and unique an element $f(a)\in\oplus_{P}
B'(P)_{p-1}$ such that $d(a+f(a))\in H'\oplus B'$.
\end{lemma}
\proof
Let us prove the uniqueness first. 
By subtraction we reduce the problem to proving that there is no nonzero $b\in
B'_{p-1}$ with $d(b)\in H'\oplus B'$. Suppose such an element $b$ exists and let
$P$ be a maximal element in $D$ with the property $b_P=\pi_P(b)\not=0$. By the
maximality of $P$ we have $d(P)(b_P)=\pi_Pd(b)\in H'(P)\oplus B'(P)$ whence 
$d(P)(b_P)=0$. Since $b_P\in B'(P)$ we conclude that $b_P=0$ which is a
contradiction.

Now we prove the existence. Using downward induction on $D$, it suffices to prove
the following; let $c\in K_p$ be such that $\pi_Pd(c)\in H'(P)\oplus B'(P)$ 
for all $P$ greater than a given $R\in D$. Then there exists $b\in B'(R)_{p-1}$ such
that $\pi_Pd(c+b)\in H'(P)\oplus B'(P)$ for $P\geq R$. 

This claim is immediate. Indeed one can take $b\in B'(R)_{p-1}$ with the
condition that $d(R)(b)=-[\pi_Rd(c)]_{B(R)}$ where $[\ \ ]_{B(R)}$ means the
projection of $K(R)$ to $B(R)$.              \qed

\medskip
Lemma \ref{map} defines a degree -1 linear (by uniquenss) map $f:H'\to B'$.
Notice that by construction if $a\in H'(Q)$ then $f(a)\in\oplus_{P<Q}B'(P)$.
We put $M=\{a+f(a)|a\in H'\}$. Clearly $M$ is a graded linear subspace of
$K$.

\begin{lemma}
\label{subcomplex}
The subspace $M$ is a subcomplex of $K$.
\end{lemma}

\proof
We need to prove that $d(M)\subset M$. 
Suppose $a\in H'$ and consider $c=d(a+f(a))$. By construction $c\in H'\oplus
B'$, i.e., $c=e+b$ where $e\in H'$ and $b\in B'$.
Since $d(c)=0\in H'\oplus B'$ we have by the uniqueness part of Lemma \ref{map}
that $b=f(e)$, i.e., $c\in M$.                      \qed

\medskip
Now we define the graded free $S$-submodule $\tilde M$ of $\tilde K$ as generated
by $M$. Lemma \ref{subcomplex} implies that $\tilde M$ is a subcomplex of $\tilde
K$. The following result is the main one of this section.

\begin{theorem}
\label{minimal}
The complex $\tilde M$ is a minimal resolution of $S/I$.
\end{theorem}

\proof
To prove this theorem it suffices to prove that the complex 
$\tilde K/\tilde M$
is exact. For that, in turn, it suffices to prove that $K'=
(\tilde K/\tilde M)\otimes k=(\tilde K\otimes k)/(\tilde M\otimes k)$
 is exact.

To analyze the complex $K'$ notice first that as graded linear spaces $\tilde K
\otimes k=K$ and $\tilde M\otimes k=M$. Moreover in the decomposition $K=B \oplus
H'\oplus B'$ we have $M\subset H'\oplus B'$ where $M$ is the graph of $f:H'\to
B'$. Thus up to natural isomorphism
$K'=B\oplus B'$ as graded linear space.
In particular 
\begin{equation}
\label{factor}
K'=\oplus_Q(B(Q)\oplus B'(Q)) 
\end{equation}
(again as graded linear spaces).
Moreover, unlike for $K$, the complexes $K(Q)$
are subcomplexes of $\tilde K\otimes k$ whence \ref{factor} holds in the category
of chain complexes, i.e., the differential in $K'$ coincides with $\oplus_Qd(Q)$. 
Now the statement follows immediately form the isomorphisms $d(Q):B'(Q)_p\to
B(Q)_{p-1}$.                     \qed

\medskip
Since in general the minimal resolution $\tilde M$ of $S/I$ is not constructed
canonically, it is
interesting to consider a case when it is canonical. This resolution was
discovered in \cite{BPS} for so called generic monomial ideals
 (see below). We want to
show how $\tilde M$ reduces to this resolution for a significantly wider class
of $\A$ (even among monomial ideals).

In \cite{BPS}, the Scarf
 complex is the subcomplex of $\tilde M$ generated by $\sigma\subset\A$
such that $|\B(Q)|=1$ for $Q=Q_{\sigma}$. 

\begin{proposition}
\label{BPS}
The complex
 $\tilde M$ coincides with the Scurf complex (in particular the latter is
acyclic) if and only if
for every $Q\in D$ either $|\B(Q)|=1$ or the complex $K(Q)$ is exact.
\end{proposition}
\proof
The condition is obviously necessary. Let us prove that it is sufficient.
Put $\bar D=\{Q\in D| |\B(Q)|=1\}$. We have $K(Q)=K_1(Q)=H'_1(Q)$ for every
$Q\in\bar D$ and 
$H'(Q)=0$ for other $Q$. Thus the following claim suffices for the proposition.

{\bf Claim}. Let $\sigma\subset \A$ be such that $Q_{\sigma}\in \bar D$. Then
$Q_{\sigma_i}\in\bar D$ for every $Q_i\in\sigma$. 

{\bf Proof of Claim.} Suppose $Q_{\sigma_i}\not\in \bar D$ for some $i$. 
Then
there exists $Q_j\in\A$ such that $Q_j$ divides $Q_{\sigma_i}$. 
There are two possibilities. One possibility is that
$Q_j\in\sigma_i$.
This implies
$Q_{\sigma_j}=Q_{\sigma}$ which is a contradiction.
The other possibility is that $Q_j\not\in\sigma_i$. If $Q_j\in\sigma$ then
$j=i$ and
$Q_{\sigma}=Q_{\sigma_i}$ which is a contradiction. If  
$Q_j\not\in\sigma$ then $Q_{\sigma}=Q_{\sigma\cup\{Q_j\}}$ which is again a
contradiction. This proves the claim and the proposition.                  \qed

\medskip
In \cite{BPS}, a monomial ideal $I$ is called generic if no variable appears with
the same nonzero exponent in two distinct minimal generators of $I$. This (though
ambiguous) term can be used for an arbitrary set $\A$ of homogeneous polynomials
(with respect ot their factorizations into irreducible factors).

\begin{proposition}
\label{boolean}
Let $\A$ be generic in the above sense. Then for every $Q\in D$
the poset $\B(Q)$ is Boolean.
\end{proposition}
\proof
The key observation is that $\B(Q)$ has a unique minimal element. Indeed write
$Q=\prod_{i=1}^rZ_i^{m_i}$ where $Z_i$ are irreducibles and $m_i>0$.
Then for each $i$ there exists a unique $Q_{j_i}\in\A$ such that $Z_i$ has
exponent $m_i$ in decomposition of $Q_{j_i}$. 
Put $\sigma=\{Q_{j_i}|i=1,\ldots,r\}$.
It is easy to see that $Q_{\sigma}=Q$ and $\sigma$ is the unique minimal element
of $\B(Q)$.

Now the result follows since $\B(Q)$ always has a unique maximal element and with
any two elements of $\B$ contains every element between them.
     \qed

\medskip
\begin{corollary}
\label{scarf}
If $\A$ is generic then $\tilde M$ coincides with the Scarf complex (cf.
\cite{BPS}).
\end{corollary}

\medskip
It is easy to find examples with nongeneric $\A$ and all Boolean $\B(Q)$ (e.g.,
$\A=\{xy,xz\})$.
Moreover Propositions \ref{BPS} and \ref{boolean} give a class of $\A$ defined
by an easily checkable condition with the
Scarf resolution.
The condition is that {\bf \centerline {$Q_{\sigma}=Q_{\tau}$
implies $Q_{\sigma\cap\tau}=Q_{\sigma}$
 for every subsets $\sigma$ and $\tau$ of $\A$.}}

The following example shows that $\B(Q)$ do not have to be Boolean to satisfy the
condition of Proposition \ref{BPS}.

{\bf Example 4.1.}
Let $\A=\{x^2yz,xy^2w,x^2zw,xy^2z\}$.
One can find easily that the condition of Proposition \ref{BPS} holds
whence
the minimal resolution coincides with the Scarf complex. On the other hand
$\B(Q)$ for
$Q=x^2y^2zw$ has three minimal elements.

Now we give a couple of examples with $H'$ not being invariant with respect to the
differential (in particular nonvanishing).

{\bf Example 4.2.}
Let $\A=\{xy,xz,yu,uv\}$ enumerated in the order they are written.
There are three $\B(Q)$ that have more than one element:
$\B(xyzu)=\{\{2,3\},\{1,2,3\}\}, \B(xyuv)=\{\{3,4\},\{1,3,4\}\}$, and \hfill\break
$\B(xyzuv)=\{\{1,2,4\},\{2,3,4\},\{1,2,3,4\}\}$. Only the last $K(Q)$ has
nontrivial homology, namely $\dim H_3=1$ and as a cycle representative of a nonzero
class one can take $\sigma=\{1,2,4\}$. Then one has
$d(\sigma)=\{1,2\}-\{1,4\}+\{2,4\}$ where $\{1,4\}$ is the trivial cycle in
$K(xyuv)$, namely the boundary of $\tau=-\{1,3,4\}$. Thus using our construction
we have $M_3=ka$ where $a=\sigma+\tau$ and then
$d(a)=\{1,2\}+\{2,4\}-\{1,3\}-\{3,4\}\in M_2$ where $M_2$ is generated by
$\{\{1,2\},\{1,3\},\{2,4\},\{3,4\}\}$.
In particular the Betti numbers of $S/I$ are 1,4,4,1.

\medskip
{\bf Example 4.3.} This is Avramov's example from \cite{Av}. Let
$$\A=\{x^2,xy,yz,zw,w^2\}$$ 
again with the natural linear order. There are six
 $\B(Q)$ with more than one element: $\B(x^2yz)=\{\{1,3\},\{1,2,3\}\}, \B(xyzw)=
\{\{2,4\},\{2,3,4\}\}, \B(yzw^2)=\{\{3,5\},\{3,4,5\}\}, \B(x^2yzw)=\{\{1,2,4\},
\{1,3,4\},\{1,2,3,4\}\}, \hfill\break
\B(xyzw^2)=\{\{2,3,5\},\{2,4,5\},\{2,3,4,5\}\}$
and \hfill\break
$\B(x^2yzw^2)=\{\{1,3,5\},\{1,2,3,5\},\{1,2,4,5\},\{1,3,4,5\},\{1,2,3,4,5\}\}$.
Only the last three give $K(Q)$ with nonvanishing homology. As cycle
representatives one can take $\{1,2,4\},\{2,4,5\}, \{1,2,4,5\}$. Then one gets
the generators  of $M_4$ as $\{1,2,4,5\}$, of $M_3$ as $\{1,2,4\}+\{2,3,4\},
\{2,4,5\}+\{2,3,4\}, \{1,2,5\},\{1,4,5\}$ and of $M_2$ as all the pairs of
generators except $\{1,3\},\{2,4\}$, and $\{3,5\}$ ($M_1=K_1$ and $M_0=K_0$ as
always). It is easy to see that $M$ is invariant under $d$ and the dimensions of
$M_p$ coincide with the respective Betti numbers of $S/I$ that are 1,5,7,4,1.
\bigskip
\section{Examples of ideals with Taylor resolution}
\bigskip
In this section we show examples of classes of sets $\A$ of homogeneous polynomials
such that the conditions of Theorem \ref{taylor} hold. The following result is
useful for that. For every $P\in W$ denote by $E_P$ the subset of $E$ of the
irreducible polynomials taking part in the factorization of $P$. Also put
$J(P)=J(E_P)$.

\begin{proposition}
\label{coatomic}
For every saturated $G\subset E$ and every $P\in W(G)$ we have $\tilde
H_p(D(G,<P))=0$ for $p\geq |E_P|-1$.
\end{proposition}
\proof
Suppose $G$ and $P$ are as in the statement and $P=\prod_{i=1}^rZ_i^{m_i}$ is the
factorization of $P$ into irreducible factors. Notice that $r=|E_P|$. Let now
$X$ and $Y$ be distinct maximal elements of $D(G,P)$. Then $\psi_G(X\vee Y)\geq P$.
This implies that for every $i$, $1\leq i \leq r$,
all but at most one maximal elements of $D(G,P)$ have the factor $Z_i^{\ell_i}$
with $\ell_i<m_i$ in
their factorizations. Since on the other hand each element of $D(G,P)$ has the
multiplicity of at least one of $Z_i$ smaller than $m_i$ we have
\begin{equation}
\label{3.1}
|\max D(G,P)|\leq r.
\end{equation}

Now the poset $D(G,P)$ can be viewed as a lattice with the maximal and minimal
elements deleted. Thus by \cite{Fo} its homology can be computed as the homology of
its coatomic complex of the lattice. The vertices of this complex are the maximal
elements of $D(G,P)$ whence \eqref{3.1} implies that either this complex is a simplex
or its dimension is less than $r-1$. The result follows.                 \qed

\medskip
\begin{corollary}
\label{regular}
If the set $E_P$ forms a regular sequence (for the module $S$)
 then the condition of Theorem
\ref{taylor} holds for $P$ and any $G$ such that $P\in \psi_G(D)$.
\end{corollary}
\proof
We have in this case $\depth J(G)\geq \depth J(P)=|E_P|.$ Thus the result follows
from the previous proposition.                   \qed

\medskip
In particular Corollary \ref{regular} recovers the well-known result that for any
monomial ideal the Taylor complex is acyclic.

Another class of ideals where the condition of Theorem \ref{taylor} simplifies
(although stays nontrivial) consists of ideals generated by products of linear
polynomials. For this class the condition of Theorem \ref{taylor} specializes to
\begin{equation}
\label{3.2}
\tilde H_p(D(G,P))=0\ {\rm for}\ p\geq \dim E_P-1
\end{equation}
where $\dim E_p$ is the usual dimension of the linear space generated by $E_P$.

We can find a sufficient condition for \eqref{3.2}
using a result similar
to (although more subtle than) Proposition \ref{coatomic} but for atomic complexes.
First we need a lemma.
\begin{lemma}
\label{homotopic2}
Let $G$ be a saturated subset of $E$ and $U$ a decreasing subset of the subposet
$\psi_G(D)$ of $W(G)$. Then $\psi_G$ is a homotopy equivalence 
of $\psi^{-1}_G(U)$ and $U$.
\end{lemma}
\proof
As in the proof of Lemma \ref{homotopic} it is easy to see that for every $P\in
U$ the poset $\psi_G^{-1}(U_{\leq P})$ has the unique maximal element
lcm$\{Q_i\in\A|\psi_G(Q_i)\leq P\}$. Thus this poset is contractable and the
result follows.                      \qed

\medskip
For each $i=1,2,\dots,m$ put $E_i=E_{Q_i}$.

\begin{proposition}
\label{atomic}
Suppose all $Q_i$ are square free products of linear polynomials. 
If for every $i=1,2,\ldots, m$ the set $G_i=E\setminus E_i$
 is saturated (that in this case means
no element of $E_i$ is a linear combination of elements from $G_i$) then \eqref{3.2}
holds.
\end{proposition}
\proof
Fix a saturated set $G\subset E$ and $P\in W(G)$. Lemma \ref{homotopic2} implies
that $\tilde H_p(D(G,P))\approx \tilde H_p(\psi_G(D)_{<P})$ for every $p$.
Let $s$ be the number of the
minimal elements $Z_{i_1},\ldots,Z_{i_s}$ of $\psi_G(D)_{<P}$, i.e., 
$Z_{i_j}=\psi_G(Q_{i_j})$ and $Z_{i_j}<P$.
If $s\leq\dim E_P$ 
 then  
the dimension of the atomic complex of $\psi_G(D)_{<P}$ is less than $\dim E_P-1$
unless this complex is a simplex (cf. the proof of
Proposition \ref{coatomic}). 
Then \eqref{3.2} follows.

Thus we have to prove only the impossibility of $s>\dim E_P$.
Suppose it is the case. Let $r$ be the largest dimension
 of $G_i$. Then the conditions
on $G_i$ and $Q_i$ imply that $\dim\bigcap_{j=1}^sG_{i_j}\leq r-s+1<r-\dim E_P+1$.
The condition $Z_{i_j}<P$ is equivalent to $G_{i_j}\supset E\setminus E_P$
whence $\bigcap_{j=1}^s G_{i_j}\supset E\setminus E_P$. We obtain
$$\dim (E\setminus E_P)< r-\dim E_P+1$$
whence
$$\dim E\leq \dim E_P+\dim(E\setminus E_P)<r+1$$
which is a contradiction.                       \qed

\medskip
The condition \eqref{3.2} simplifies significantly if we assume that the set $E$ of
linear polynomials is generic, i.e., any subset of $E$ with at most $n$ elements
is linearly independent. 
\begin{proposition}
\label{generic}
Suppose that all $Q_i$ are products of linear polynomials and $E$ is generic.
Then \eqref{3.2} is equivalent to
\begin{equation}
\label{3.2'} 
\oplus_{P\in D}\tilde H_p(D_{<P})=0\ {\rm for}\ p\geq n-1.
\end{equation}
Equality \eqref{3.2'} for $P\in D$ with $\rk
E_P=n$ suffices for \eqref{3.2}.
\end{proposition}
\proof
Suppose first that $G$ is a proper saturated subset of $E$ and $P\in W(G)$.
Since $E$ is generic, $G$ is linearly independent and so is $E_P\subset G$,
 i.e., $\rk E_P=|E_P|$. Proposition \ref{coatomic} implies \eqref{3.2}
for this $G$ and $P$.

Now suppose that $G=E$ whence $\psi_G$ is the canonical embedding $D\subset W$.
Then the only nontrivial case is where $P\in D$ whence $D(G,<P)=D_{<P}$.
 If $\rk E_P<n$ then $E_P$ is
again linearly independent and \eqref{3.2} holds by the same reason as in the
previous paragraph. If $\rk E_P=n$ then \eqref{3.2} is obviously equivalent to
\eqref{3.2'}.                              \qed

\medskip
Notice that equality \eqref{3.2'} using only the combinatorics of the lattice $D$
and $n$. This combinatorics can be expressed in terms of
 other polynomial ideals, the
simplest from them being monomial ones. More precisely assign to each linear
polynomial $Z_i\in E$ ($i=1,2,\ldots,r$) an indeterminant $y_i$ and consider the
polynomial ring $\tilde S=k[y_1,\ldots,y_r]$. The natural algebra
map $\tilde
S\to S$ via $y_i\mapsto Z_i$ assigns to each $Q_j\in \A$ a monomial $\tilde
Q_j\in\tilde S$ which generate the monomial ideal $\tilde I$ of $\tilde S$.
Then the following result follows straightforwardly from the results of the
previous sections.

\begin{corollary}
\label{monoms}
Undeer the conditions of Proposition \ref{generic} the Taylor complex of $S/I$ is a
resolution if and only if $b_p(\tilde S/\tilde I)=0$ for $p\geq n+1$. If this
condition holds then the $k$-algebras ${\rm Tor}_*^S(S/I,k)$ and ${\rm
Tor}_*^{\tilde S}(\tilde S/\tilde I,k)$ are naturally isomorphic.
\end{corollary}

\medskip
Using Proposition \ref{generic}, 
we can easily give a series of examples of non-\hfill\break
monomial ideals whose Taylor
complexes are exact.

{\bf Example 5.1.}
 Consider generic arrangement of more than $n$ hyperplanes in a space of
dimension $n$.
Let $\A$ consist of no more than $n$ arbitrary
products of functionals of hyperplanes.
Then the Taylor complex of $\A$ is its resolution.
Indeed for every $P\in D_0$ the poset $D_{<P}$ has at most $n$ atoms whence
its atomic complex is either a simplex or has dimension less than $n-1$. The
condition \eqref{3.2'} follows.

\medskip
\begin{remark}
In the case where $Q_i$ are square free products of linear polynomials
the atomic complex of $D_{<P}$ can be interpreted as the nerve of
the collection $\{G_i=E_P\setminus E_i|E_i\subset E_P\}$ of subsets of $E_P$.
This is the reason for the relations between the Betti numbers of
$S/I$ and the complement of a coordinate subspace arrangement studied for
monomial ideals in \cite{GPW1}. In particular the duality between $E_i$ and $G_i$
leads to the appearance of the
Alexander dual complexes.
\end{remark}

\end{document}